\documentclass[a4paper,20pt]{amsproc}

\usepackage{bbm}
\usepackage{hyperref}
\newtheorem{theorem}{Theorem}[section]

\newtheorem{corollary}[theorem]{Corollary}

\theoremstyle{definition}
\newtheorem{definition}[theorem]{Definition}

\theoremstyle{remark}

\numberwithin{equation}{section}

\begin{document}
\LARGE
\title{On a conjecture of the paper Differentiability of the Conjugacy in the Hartman-Grobman Theorem}


\author[Genrich Belitskii]{Genrich Belitskii}
\email{genrich@cs.bgu.ac.il}

\author[Victoria Rayskin]{Victoria Rayskin}
\email{vrayskin@gmail.com}


\date{}

\begin{abstract}
In this note we show that for the construction of differentiable conjugation, the assumption of the existence of smooth bump function is not necessary, and consequently the corresponding conjecture stated in the paper of W. Zhang, K. Lu and W. Zhang ``Dfifferentiability of the Conjugacy in the Hartman-Grobman Theorem'' (\cite{ZLZ}) is incorrect. We show that instead of bump functions we can use smooth blid maps. We also propose a construction of the blid map for the space $X=C^0[0,1]$, which does not possess a smooth bump function.

\end{abstract}

\maketitle

\section{Introduction}\label{sec-intro}

Linearization and normal forms are convenient simplification of complex dynamics. 

For a diffeomorphism $F$ with a fixed point $0$, we would like to find a smooth transformation $\Phi$ defined in a neighborhood of $0$ such that $\Phi\circ F \circ \Phi^{-1}$ has a simplified (polynomial) form called the normal form. If $\Phi\circ F \circ \Phi^{-1}=DF=\Lambda $ is linear, the conjugation is called linearization. There are two major questions in this area of research: how to increase smoothness of the conjugation $\Phi$, and whether it is sufficient to assume low smoothness of the diffeomorphism $F$.

Hartman and Grobman independently showed that if $\Lambda$ is hyperbolic, then for a diffeomorphism $F$ there exists a local homeomorphism $\Phi$ such that $\Phi \circ F \circ \Phi^{-1}=\Lambda$. Different proofs were given by Pugh \cite{P}. A higher regularity of $\Phi$ has been an active area of research.

The first attempt to answer the question of differentiability of $\Phi$ at the fixed point $0$ under hyperbolicity assumption was made in [48], but an error was found and discussed in \cite{R}. Later, in \cite{GHR}, Guysinsky, Hasselblatt and Rayskin presented correct proof. However, it was restricted to $F\in C^{\infty}$ (or more precisely, it was restricted to $F\in C^k$, where $k$ is defined by complicated expression). It was conjectured in the paper that the result is correct for $F\in C^2$, as it was announced in \cite{vS}.

In the paper of Zhang, Lu and Zhang (\cite{ZLZ}) it was shown that for a Banach space diffeomorphism $F$ with a hyperbolic fixed point and $\alpha$-H{\"o}lder $DF$, the local conjugating homeomorphism $\Phi$
is differentiable at the fixed point. Moreover, 
$$
\Phi(x) = x+O(||x||^{1+\beta}) \mbox{\  and \ } \Phi^{-1}(x) = x+O(||x||^{1+\beta})
$$
as $x\to 0$, for certain $\beta \in (0,\alpha]$. 

There are two additional assumptions in this theorem. The first one is the spectral band width inequality. It is explained that this inequality is sharp if the spectrum has at most one connected component inside of the unit circle in $X$, and at most one connected component outside of the unit circle in $X$.
%
%
%
%
It is pointed out in \cite{ZLZ} that this is not a non-resonance condition required for smooth linearization.

The second assumption is the assumption that the Banach space must possess smooth bump functions. 
It is conjectured in the paper that the second assumption is a necessary condition. 

Below we explain that this conjecture is not correct. The bump function condition can be replaced with blid map condition, which is less restrictive as explained below.

The Banach-Mazur theorem states that any real separable Banach space is isometrically isomorphic to a closed subspace of $C[0,1]$.  Consequently, the space of $C[0,1]$ does not have smooth bump functions at all 
(see~\cite{K} or \cite{M}). 

For localization of some Banach spaces and sets that do not possess bump functions, in our paper \cite{BR2}, we introduce blid maps:
\begin{definition}\label{def-blid map} A  $C^q$-{\it blid map} for a Banach space $X$ is a global {\bf B}ounded {\bf L}ocal {\bf Id}entity at zero $C^q$-map $H:X \to X$.
\end{definition}

     In other words, $H$ is a global representative of the germ at zero of identity map such that
 $$
                               \sup_x || H(x)||<\infty.
 $$
Blid maps allow locally defined mappings to be extended to the whole space. The idea of using bounded local identity maps was first discussed in \cite{B}. Later, blid maps were employed in \cite{BR3}. 
 We will show how to use blid maps for the proof of Theorem 7.1 of \cite{ZLZ}. Application of blid maps to this theorem extends the result to a bigger variety of Banach spaces, because blid maps exist on any space with bump function (multiplication by Identity, turns bump function into a blid map). Blid maps allow to reformulate Theorem 7.1 in the following way:
\begin{theorem}\label{thm-diff}
Let $X$ be a Banach space possessing a $C^1$-blid map with bounded derivative.
Suppose $F:X\to X$ is a diffeomorphism with a hyperbolic fixed point,  $DF$ is $\alpha$-H{\"o}lder, and the spectral band width condition is satisfied.
Then, the local conjugating homeomorphism $\Phi$
is differentiable at the fixed point. Moreover, 
$$
\Phi(x) = x+O(||x||^{1+\beta}) \mbox{\  and \ } \Phi^{-1}(x) = x+O(||x||^{1+\beta})
$$
as $x\to 0$, for certain $\beta \in (0,\alpha]$. 
\end{theorem}
Now we will also show how to construct blid maps for $C[0,1]$. For other examples of construction of blid maps we refer the reader to \cite{BR2}.

Let  $X=C[0,1]$ be the space of all continuous functions on $[0,1]$ with
 $$
                     ||x||=\max_t|x(t)| ,\ \  t\in [0,1],
 $$
and let h be a $C^\infty$-bump function on the real line. Then the map
 
 \begin{equation}\label{C[01]-blid map}
                                     H(x)(t)=h(x(t))x(t),\ \  x\in X
 \end{equation}
is $C^\infty$-blid map with bounded derivatives of all orders. 

Indeed, since $h$ is locally equal to 1, $H$ is a local identity. Let $\epsilon$ be a positive real number such that $h(\tau)\equiv 0$ for all $|\tau| >\epsilon$. We can always find such $\epsilon$, because bump functions have bounded support.\\ 
Then, $$||H(x)(t)|| = || h(x(t))x(t) ||\leq \epsilon.$$ Also, \\
$$
||H'(x)(t)||= ||h'(x(t)) x(t)+ h(x(t))|| \leq \epsilon \sup|h'| +1.
$$
\\
Similarly, one can show boundedness of all higher order derivatives.
\begin{corollary}
Let $X=C[0,1]$.
Suppose $F:X\to X$ is a diffeomorphism with a hyperbolic fixed point,  $DF$ is $\alpha$-H{\"o}lder, and the spectral band width condition is satisfied.
Then, the local conjugating homeomorphism $\Phi$
is differentiable at the fixed point. Moreover, 
$$
\Phi(x) = x+O(||x||^{1+\beta}) \mbox{\  and \ } \Phi^{-1}(x) = x+O(||x||^{1+\beta})
$$
as $x\to 0$, for certain $\beta \in (0,\alpha]$. 
\end{corollary}

\section{Proof of Theorem \ref{thm-diff}}
Zhang, Lu and Zhang showed that for the conclusion of their Theorem 7.1 it is enough to satisfy the inequalities 1 and 2 (see \ref{ineq} below),  which are called condition (7.6) in their paper.

In order to apply the blid maps instead of bump functions to the inequalities (\ref{ineq}), it is sufficient to construct a bounded blid map, which has only first-order bounded derivative. I.e., let blid map $H(x): X\to X$ be as follows:
\begin{equation*}
\begin{array}{l}
\mbox{1. \ }  H(x) = x \mbox{\ for\ } ||x||<1\\
\mbox{2. \ } H\in C^1 \mbox{\ and\ } ||H^{(j)}(x)||\leq c_j,  \ j=0,1 .
\end{array}
\end{equation*} 

The condition (7.6) of \cite{ZLZ} is:
\begin{equation}\label{ineq}
\begin{array}{l}
\mbox{1. \ } \sup_{x\in X}||DF(x) -\Lambda|| \leq \delta_{\eta}\\
\mbox{2. \ } \sup_{x\in V\setminus O}\left\{||DF(x) - \Lambda|| / ||x||^{\alpha}\right\} = M < \infty 
\end{array}
\end{equation} 

Let $DF -\Lambda = f$. Define
$$
\tilde{f}(x):= f\left( \delta H(x/ \delta) \right)
$$

We will show that is $f$ satisfies (7.6), then so does $\tilde{f}$.
$$
\sup_{x\in X}|| D\tilde{f}(x)|| \leq \sup_{x\in X} || D f(x) || \cdot \sup_{x\in X} || DH(x) || \leq \delta_{\eta} \cdot c_1.
$$
Thus, the first inequality of (7.6) holds for $\tilde{f}$. For the second inequality we have the following estimate:
$$
\frac{|| D \tilde{f}(x)||}{||x||^{\alpha}} \leq \frac{|| Df \left(\delta H(x/\delta)\right) ||}{|| \delta H(x/\delta) ||^{\alpha}} \cdot \left( \frac{|| \delta H(x/\delta) ||}{||x||}  \right)^{\alpha}.
$$
The second multiple is bounded, because for small $x$ (say, $||x/\delta||<\epsilon$ for some $\epsilon>0$) we have
$$
\frac{|| \delta H(x/\delta) ||}{||x||} < c_1 + o(1),
$$
while for $||x/\delta|| \geq \epsilon$
$$
\frac{|| \delta H(x/\delta) ||}{||x||} < c_0/\epsilon.
$$
I.e., $\frac{|| \delta H(x/\delta) ||}{||x||}$ is less than some constant $m$.
Then, 
$$ 
\sup_{x\in V\setminus O}\frac{|| D \tilde{f}(x)||}{||x||^{\alpha}} \leq$$ $$
\sup_{0<||x||<\delta c_0}\left\{||D f(x)|| / ||x||^{\alpha}\right\} \cdot \sup_{x\in X}|| D H(x)|| \cdot m^{\alpha}=$$ $$ \sup_{0<||x||<\delta c_0}\left\{||D f(x)|| / ||x||^{\alpha}\right\}c_1\cdot m^{\alpha}
$$

This quantity is bounded by $M c_1 m^{\alpha}$ if $\delta$ is sufficiently small.
\bibliographystyle{amsalpha}

\end{document}